\documentclass[number,citesort,dvips]{arxbj}

\aid{0}
\volume{17}
\issue{1}
\pubyear{2011}
\firstpage{424}
\lastpage{440}
\doi{10.3150/10-BEJ276}

\makeatletter
\newcommand{\eqref}[1]{(\ref{#1})}
\newtheorem{theorem}{Theorem}
\newtheorem{proposition}[theorem]{Proposition}
\newtheorem{corollary}[theorem]{Corollary}

\newremark{remark}[theorem]{Remark}
\makeatother

\begin{document}
\begin{frontmatter}

\title{On the functional central limit theorem via martingale approximation}
\runtitle{Functional CLT via martingale approximation}

\begin{aug}
\author[a]{\fnms{Mikhail} \snm{Gordin}\thanksref{a}\ead[label=e1]{gordin@pdmi.ras.ru}\corref{}}
\and
\author[b]{\fnms{Magda} \snm{Peligrad}\thanksref{b}\ead[label=e2]{peligrm@ucmail.uc.edu}}
\runauthor{M. Gordin and M. Peligrad}
\address[a]{POMI (Saint Petersburg Department of the Steklov Institute of
Mathematics), 27 Fontanka emb., Saint Petersburg 191023, Russia.
\printead{e1}}
\address[b]{Department of Mathematical Sciences, University of
Cincinnati, P.O. Box
210025, Cincinnati, OH 45221-0025, USA. \printead{e2}}
\end{aug}

% HISTORY:
\received{\smonth{7} \syear{2009}}
\revised{\smonth{2} \syear{2010}}

% ABSTRACT
%
\begin{abstract}
In this paper, we develop necessary and sufficient conditions for the
validity of a martingale approximation for the partial sums of a stationary
process in terms of the maximum of consecutive errors. Such an approximation
is useful for transferring the
conditional functional central limit theorem from the martingale to the
original process. The condition found is simple
and well adapted to a variety of examples, leading to a better understanding
of the structure of several stochastic processes and their asymptotic
behaviors. The approximation brings together many disparate examples in
probability theory. It is valid for classes of variables defined by familiar
projection conditions such as the Maxwell--Woodroofe condition, various classes
of mixing processes, including the large class of strongly mixing processes,
and for additive functionals of Markov chains with normal or symmetric
Markov operators.

\end{abstract}

% KEYWORDS
%

\end{frontmatter}
%

%s1 ###
\section{Introduction and results}

The objective of this paper is to find a characterization of stationary
stochastic processes that can be studied via a martingale approximation in
order to derive the functional central limit theorem for processes
associated with partial sums.

There are several ways to present the results since stationary processes
can be introduced in several equivalent ways. We assume that $(\xi
_{n})_{n\in%
\mathbb{Z}}$ denotes a stationary Markov chain defined on a probability
space $(\Omega,\mathcal{F},P)$ with values in a measurable space $(S,%
\mathcal{A})$. The marginal distribution and the transition kernel are
denoted by $\pi(A)=P(\xi_{0}\in A)$ and $Q(\xi_{0},A)=P(\xi_{1}\in A|
\xi
_{0})$, respectively. In addition, $Q$ denotes the operator {acting via
$(Qf)(\xi)=\int
_{S}f(s)Q(\xi,\mathrm{d}s).$ Next, let $\mathbb{L}_{0}^{2}(\pi)$ be the set of
functions on $S$ such that $\int f^{2}\,\mathrm{d}\pi<\infty$ and $\int f\,\mathrm{d}\pi=0.$
Denote by $\mathcal{F}_{k}$ the $\sigma$-field generated by $\xi_{i}$
with $%
i\leq k,$ $X_{i}=f(\xi_{i})$, $S_{n}=\sum _{i=0}^{n-1}X_{i}$
(i.e., $%
S_{0}=0,S_{1}=X_{0}$, $S_{2}=X_{0}+X_{1},\dots$). For any integrable variable
$X$, we define $\mathbb{E}_{k}(X)=\mathbb{E}(X|\mathcal{F}_{k}).$ In our
notation, $\mathbb{E}_{0}(X_{1})=Qf(\xi_{0})=\mathbb{E}(X_{1}|\xi
_{0}).$ }
We also set $\mathcal F_{-\infty}=\bigcap_{k \in\mathbb Z}
\mathcal F_{k} $ .

Throughout the paper, we assume $f\in\mathbb{L}_{0}^{2}(\pi)$; in other
words, we assume that $\Vert X\Vert _{2}=(\mathbb{E}[X_{1}^{2}])^{1/2}<\infty$
and $%
\mathbb{E}[X_{1}]=0.$

Note that any stationary sequence $(Y_{k})_{k\in\mathbb{Z}}$ can be viewed
as a function of a Markov process $\xi_{k}=(Y_{i};i\leq k)$ for the
function $g(\xi_{k})=Y_{k}$.

The stationary stochastic processes may be also introduced in the following,
alternative, way. Let $T\dvtx \Omega\mapsto\Omega$ be a bijective bimeasurable
transformation preserving the probability. Let $\mathcal{F}_{0}$ be a sub-$%
\sigma$-algebra of $\mathcal{F}$ satisfying $\mathcal{F}_{0}\subseteq
T^{-1}(\mathcal{F}_{0})$. We then define the non-decreasing filtration
$(%
\mathcal{F}_{i})_{i\in\mathbb{Z}}$ by $\mathcal{F}_{i}=T^{-i}(\mathcal
{F}%
_{0})$. Let $X_{0}$ be a random variable which is ${}\mathcal{F}_{0}$%
-measurable. We define the stationary sequence $(X_{i})_{i\in\mathbb
{Z}}$ by
$X_{i}=X_{0}\circ T^{i}$.

In this paper, we shall use both frameworks.

In order to analyze the asymptotic behavior of the partial sums $%
S_{n}=\sum_{i=0}^{n-1}X_{i},$ Gordin, in~\cite{g}, proposed to
decompose the
sums related to the original stationary sequence into the sum
%
%e1 ###
\begin{equation}
S_{n}=M_{n}+R_{n} \label{martaprox}
\end{equation}
of a square-integrable martingale $M_{n}=\sum_{i=0}^{n-1}D_{i}$ adapted
to $%
\mathcal{F}_{n}$, whose martingale differences~$(D_{i})$ are
stationary, and
a so-called coboundary $R_{n}$, that is, a telescoping sum of random variables
with the basic property that $\sup_{n}\mathbb{E}(R_{n}^{2})<\infty$. More
precisely, $X_{n}=D_{n}+Z_{n}-Z_{n-1},$ where $Z_{n}$ is another stationary
sequence in $\mathbb{L}_{2}$. The limiting properties of the
martingales can then be transported from the martingale to the general
sequence. In
the context of Markov chains, the existence of such a decomposition is
equivalent to the solvability of the Poisson equation in~$\mathbb{L}_{2}$.

  For proving a central limit theorem for stationary sequences, a
weaker form of martingale
approximation has been pointed out by many authors (see, e.g., \cite{mpu2} for a survey). Recently, two interesting
papers, one by Dedecker, Merlev\`{e}de and Voln\'{y} \cite{dmv} and
the other by
Zhao and Woodroofe \cite{Zw2}, provided necessary and sufficient
conditions for
martingale approximation with an error term in (\ref{martaprox}) satisfying
%
%e2 ###
\begin{equation}
\mathbb{E}\bigl((S_{n}-M_{n})^{2}\bigr)/n\rightarrow0. \label{MA}
\end{equation}
This decomposition is strong enough for transporting the conditional central
limit theorem from sums of stationary martingale differences in $%
\mathbb{L}_{2}$ to $S_{n}/\sqrt{n}.$ By conditional CLT, as
discussed in~\cite{DM02}, we
understand, in this context, that for any continuous function $f$ such that
$|f(x)|/(1+x^{2})$ is bounded and for any $k \ge0,$
%
%e3 ###
\begin{equation}
\biggl\Vert \mathbb{E}_{k}\bigl(f\bigl(S_{n}/\sqrt{n}\bigr)\bigr)- \int_{- \infty}^{\infty} f\bigl(x
\sqrt{\eta}\bigr) g(x)\,\mathrm{d}x\biggr\Vert _1 \mathop{\longrightarrow}_{n \to\infty} 0,
\label{CLT}
\end{equation}
where $g$ is the standard normal density and $\eta\ge0$ is an
invariant function satisfying
\[
\mathop{\lim}_{n \to\infty}\biggl\Vert\frac{\mathbb
{E}_{0}(S_n^2)}{n}-\eta
\biggr\Vert_1=0.
\]
Here, and throughout the paper, we
denote by $\Vert \cdot\Vert _{p}$ the norm in $\mathbb{L}_{p}$.

An important extension of this theory is to consider the conditional central
limit theorem in its functional form. For $t\in\lbrack0,1]$, define
\[
S_n(t)=S_{[nt]}+(nt-[nt])X_{[nt]},
\]
where $[x]$ denotes the integer part of $x$. Note that
$S_n(\cdot)/\sqrt{n}$ is a random element of the space $C([0,1])$
endowed with the supremum norm $\Vert \cdot\Vert _{\infty}.$
%the space of all functions on $[0,1]$ which have left-hand limits and
%are continuous from the right.
Then, by the conditional CLT in the functional form (FCLT), we
understand that
for any continuous function $f\dvtx  C([0,1]) \to\mathbb R$ such that
$x \mapsto|f(x)|/(1+\Vert x\Vert _{\infty}^2)$ is bounded and for any $k\ge
0$,
we have
%
%e4 ###
\begin{equation}
\biggl\Vert \mathbb{E}_{k}\bigl(f\bigl(S_n/\sqrt{n}\bigr)\bigr)- \int_{C([0,1])}\bigl(f\bigl(x
\sqrt{\eta}\bigr)\bigr)\,\mathrm{d}W(x)\biggr\Vert _1\mathop{\longrightarrow}_{n \to\infty} 0. \label{FCLT}
\end{equation}
Here, $W$ is the standard Wiener measure on $C([0,1])$.

It is well known that a martingale with stationary differences in
$\mathbb{L}_{2}$ satisfies this type of behavior with $\eta=
\lim_{n \to\infty} \sum_{l=0}^{n-1}D_{l}^2/n$ in $\mathbb
{L}_{1}$ -- this is at the heart of many statistical
procedures. This conditional form of the invariance principle is a
stable type of convergence that makes possible the change of
measure with another absolutely continuous measure, as discussed
in \cite{b,ru,hh}.

With such a result in mind, the question is now to find necessary and
sufficient conditions for a martingale decomposition with the error term
satisfying
%
%e5 ###
\begin{equation}
\mathbb{E}\Bigl(\max_{1\leq j\leq n}(S_{j}-M_{j})^{2}\Bigr)\big/n\rightarrow0 .
\label{maxcond}
\end{equation}

In order to state our martingale approximation result, for fixed $m$, we
consider the stationary sequence
%
%e6 ###
\begin{equation}
Y_{0}^{m}=\frac{1}{m}\mathbb{E}_{0}(X_{1}+\cdots+X_{m}),\qquad
Y_{k}^{m}=Y_{0}^{m}\circ T^{k}. \label{defY}
\end{equation}
In the language of Markov operators, we then have
\[
Y_{0}^{m}=\frac{1}{m}(Qf+\cdots+Q^{m}f)(\xi_{0}) .
\]
It is convenient to introduce a seminorm notation, namely,%
\[
\Vert Z\Vert _{M^{+}}=\lim\sup_{n\rightarrow\infty}\frac{1}{\sqrt{n}}\Biggl\Vert \max_{1\leq k\leq n}\Biggl| \sum_{j=1}^{k} Z \circ T^j
\Biggr|  \Biggr\Vert_{2}
\]
on the space of all $Z \in L^2_0$ with $\Vert Z\Vert _{M^{+}} < \infty.$
\begin{theorem}
\label{T}\label{maxgencopy(2)}Assume that $(X_{k})_{k\in Z}$ is a stationary
sequence of centered square-integrable random variables. Then%
%
%e7 ###
\begin{equation}
\Vert Y_{0}^{m}\Vert _{M^{+}}\rightarrow0\qquad\mbox{as }m\rightarrow\infty
 \label{MAX}
\end{equation}
if and only if there exists a martingale with stationary increments
satisfying \eqref{maxcond}. Such a martingale is unique if it exists.
In particular,
\eqref{MAX} implies \eqref{FCLT}.
\end{theorem}

As a consequence of the proof of Theorem \ref{T}, we also obtain the
following result that adds a new equivalent condition to the
characterizations by Dedecker, Merlev\`{e}de and Voln\'{y} \cite{dmv} and
Zhao and Woodroofe \cite{Zw2}. With $(Y_{k}^{m})_{k\in Z}$ defined by %
\eqref{defY} and the seminorm notation
\[
\Vert Y_{0}^{m}\Vert _{+}=\lim\sup_{n\rightarrow\infty}\frac{1}{\sqrt{n}}%
\Biggl\Vert \sum_{j=1}^{n}Y_{j}^{m}\Biggr\Vert _{2}
\]
we have the following characterization.

\begin{theorem}
\label{T2}\label{maxgencopy(1)}\label{Tcopy(1)}\label{maxgen}Assume
that $%
(X_{k})_{k\in Z}$ \ is as in Theorem \ref{T}. Then
%
%e8 ###
\begin{equation}
\Vert Y_{0}^{m}\Vert _{+}\rightarrow0\qquad\mbox{as }m\rightarrow\infty
\label{MA1}
\end{equation}
if and only if there exists a stationary martingale satisfying \eqref{MA}.
Such a martingale is unique if it exists. In particular,
\eqref{MA1} implies \eqref{CLT}.
\end{theorem}

Our approach is constructive. If the stationary sequence is supposed to be
ergodic, then the constructed martingale differences are also ergodic and
therefore the conditional theorems (\ref{CLT}) and (\ref{FCLT}) can be easily
transported to the original processes satisfying (\ref{MA1}) and (\ref
{MAX}), respectively, with $\eta=\Vert D_{0}\Vert _{2}$.

A natural and useful question is to provide classes of stochastic processes
that have a martingale decomposition with an error term satisfying %
\eqref{maxcond}, in other words, to provide sharp sufficient conditions for
such a decomposition. Obviously, a maximal inequality is needed in
order to
verify this condition. We shall combine our approach with several
maximal inequalities. One is due to Rio \cite{rio}, formula (3.9),
page 53;
for related inequalities, see \cite{Pel} and \cite{DR}.

\begin{itemize}
\item For any stationary process with centered variables in $\mathbb{L }_{2}
$,
%
%e9 ###
\begin{equation}
\mathbb{E}\Bigl(\max_{1\leq i\leq n}S_{i}^{2}\Bigr)\leq8n\mathbb{E}(X_{0}^{2})+16
\sum_{k=2}^{n}\mathbb{E}|X_{0}\mathbb{E}_{0}(S_{k}-S_1)|.
\label{Rio}
\end{equation}
\end{itemize}

Another inequality comes from \cite{PU1},
Proposition (2.3);
see also \cite{puw}, Theorem 1, for the
inequality in $%
\mathbb{L}_{p}$.

\begin{itemize}
\item For any stationary process with centered variables in $\mathbb{L}_{2}$,
\begin{eqnarray}\label{PU}
\mathbb{E}\Bigl(\max_{1\leq i\leq n}S_{i}^{2}\Bigr) & \leq& n\Biggl(2\Vert
X_{0}\Vert_{2}+3\sum_{j=0}^{r-1}\frac{\Vert \mathbb
{E}_{0}(S_{2^{j}})\Vert _{2}}{%
2^{j/2}}\Biggr)^{2}
\nonumber
\\[-8pt]
\\[-8pt]
\nonumber
& \leq& n\Biggl(2\Vert X_{0}\Vert_{2}+80\sum_{j=1}^{n}\frac{\Vert \mathbb{E}%
_{0}(S_{j})\Vert _{2}}{j^{3/2}}\Biggr)^{2},
\end{eqnarray}
where $2^{r-1} < n\leq2^{r}.$
\end{itemize}

The following maximal inequality is a particular case of Dedecker and
Merlev\`{e}de \cite{DM02}, Proposition 6; see \cite{Wu}, Theorem 1,
for the
inequality in $\mathbb{L}_{p}$.

\begin{itemize}
\item For any stationary process with centered variables in $\mathbb{L}_{2}$
such that $\mathbb E(X_{0}|\mathcal{F}_{-\infty})=0$ almost surely, we
have%
%
%e10 ###
\begin{equation}
\mathbb{E}\Bigl(\max_{1\leq i\leq n}S_{i}^{2}\Bigr)\leq4n\Biggl(\sum
_{i=0}^{\infty}\Vert
\mathbb{E}_{-i}(X_{0})-\mathbb{E}_{-i-1}(X_{0})\Vert_{2}\Biggr)^{2}
.\label{DM}
\end{equation}
\end{itemize}

Another inequality we use for additive functionals of stationary reversible
Markov chains is a consequence of Wu \cite{wu}, Corollary 2.7 and relation
(2.5) in the same paper (note that there is a typographical error in
this relation, namely, a
square should be added to the norm); see also \cite{svy}:

\begin{itemize}
\item Assume $(\xi_{n})_{n\in\mathbb{Z}}$ is a stationary,
reversible Markov chain and $X_n=f(\xi_n)$ with $f\in\mathbb
{L}_{0}^{2}(\pi)$. Then, for every $n\geq1$,%
%
%e11 ###
\begin{equation}
\mathbb{E}\Bigl(\max_{1\leq i\leq n}S_{i}^{2}\Bigr)\leq
(24n+3)\sum _{n=0}^{\infty}\mathbb E(X_{0}X_{n}),
\label{LW}
\end{equation}
provided the series on the right-hand
side is convergent.
\end{itemize}

This inequality, originally stated for the ergodic case, extends
without changes
to the general case.

By combining the martingale decomposition in Theorem \ref{T} with these
maximal inequalities, we
point out various classes of stochastic processes for which a conditional
functional limit theorem holds. These include mixing processes and
classes of
Markov chains.

%s2 ###
\section[Proof of Theorem 1]{Proof of Theorem \protect\ref{T}}

The proof of this theorem has several steps.

\textit{Step \textup{1.} Construction of the approximating martingale}.

The construction of the martingale decomposition is based on averages. It
was introduced by Wu and Woodroofe \cite{W} (see their
definition (6) on page 1677) and further developed in~\cite{Zw2}, extending the construction in \cite{Heyde} and
\cite{gl}; see also \cite{1994}, Theorem 8.1, and \cite{KV}. We give
the martingale
construction here for completeness.

We introduce a parameter $m\geq1$ (kept fixed for the moment) and define
the following stationary sequence of random variables:%
\[
\theta_{0}^{m}=\frac{1}{m}\sum_{i=1}^{m}\mathbb{E}_{0}(S_{i}),%
\theta_{k}^{m}=\theta_{0}^{m}\circ T^{k}.
\]
Set
%
%e12 ###
\begin{equation}
D_{k}^{m}=\theta_{k+1}^{m}-\mathbb{E}_{k}(\theta_{k+1}^{m}),\qquad
M_{n}^{m}=\sum_{k=0}^{n-1}D_{k}^{m}. \label{martd}
\end{equation}
Then $(D_{k}^{m})_{k\in\mathbb{Z}}$ is a stationary martingale difference
sequence and $(M_{n}^{m})_{n\geq0}$ is a martingale. Thus, we have%
\[
X_{k}=D_{k}^{m}+\theta_{k}^{m}-\theta_{k+1}^{m}+\frac{1}{m}\mathbb{E}%
_{k}(S_{k+m+1}-S_{k+1})
\]
and therefore%
\begin{eqnarray} \label{martdec}
S_{k}& =&M_{k}^{m}+\theta_{0}^{m}-\theta_{k}^{m}+\sum
_{j=1}^{k}%
\frac{1}{m}\mathbb{E}_{j-1}(S_{j+m}-S_{j})
\nonumber
\\[-8pt]
\\[-8pt]
\nonumber
& =&M_{k}^{m}+\theta_{0}^{m}-\theta_{k}^{m}+\overline{R}_{k}^{ m},
\end{eqnarray}
where we have made use of the notation%
\[
\overline{R}_{k}^{ m}=\sum _{j=1}^{k}\frac{1}{m}\mathbb{E}%
_{j-1}(S_{j+m}-S_{j}).
\]
Observe that
%
%e13 ###
\begin{equation}
\overline{R}_{k}^{ m}=\sum_{j=0}^{k-1}Y_{j}^{m}.\label{rests}
\end{equation}
With the notation
%
%e14 ###
\begin{equation}
R_{k}^{m}=\theta_{0}^{m}-\theta_{k}^{m}+\overline{R}_{k}^{ m},
\label{rest}
\end{equation}
we have
%
%e15 ###
\begin{equation}
S_{k}=M_{k}^{m}+R_{k}^{m}. \label{deco}
\end{equation}

\textit{Step \textup{2.} Sufficiency}.

We show that $\Vert Y_{0}^{m}\Vert _{M^{+}}\rightarrow0$ as $m\rightarrow\infty
$ is
sufficient for (\ref{maxcond}).

The starting point is the construction of the martingale differences,
as in (\ref{martd}). By the
martingale property and (\ref{deco}), for all positive integers
$m^{\prime}
$ and $m^{\prime\prime}$, we have
\[
\Vert D_{0}^{m^{\prime}}-D_{0}^{m^{\prime\prime}}\Vert _{2}=\frac{1}{\sqrt
{n}}%
\Vert M_{n}^{m^{\prime}}-M_{n}^{m^{\prime\prime}}\Vert _{2}=\frac{1}{\sqrt
{n}}%
\Vert R_{n}^{m^{\prime}}-R_{n}^{m^{\prime\prime}}\Vert _{2} .
\]
We now let $n\rightarrow\infty.$ By relation (\ref{rest})
and stationarity, it follows that
\begin{eqnarray*}
\lim\sup_{n\rightarrow\infty}\frac{1}{\sqrt{n}}\Vert R_{n}^{m^{\prime
}}-R_{n}^{m^{\prime\prime}}\Vert _{2}&=&\lim\sup_{n\rightarrow\infty
}\frac{1}{%
\sqrt{n}}\Vert \overline{R}_{n}^{ m^{\prime}}-\overline{R}_{n}^{
m^{\prime
\prime}}\Vert _{2} \\
&\leq&\lim\sup_{n\rightarrow\infty}\frac{1}{\sqrt{n}}(\Vert \overline{R}%
_{n}^{ m^{\prime}}\Vert _{2}+\Vert \overline{R}_{n}^{ m^{\prime\prime
}}\Vert _{2})%
.
\end{eqnarray*}
By (\ref{MAX})$,$ the limit when $m^{\prime}$ and $m^{\prime\prime}$ both
tend to $\infty$ is then $0$, giving that $(D_{0}^{m})$ is Cauchy in $%
\mathbb{L}_{2}$ and therefore convergent. Denote its limit by $D_{0}$.
Then $%
M_{n}=\sum _{k=0}^{n-1}D_{k}$ is a martingale with the desired
properties. To see this, we start from the decomposition in relation
(\ref%
{martdec}) and obtain
\[
|S_{k}-M_{k}|\leq|M_{k}^{m}-M_{k}|+|\theta_{k}^{m}-\theta_{0}^{m}|+|%
\overline{R}_{k}^{m}|.
\]
Then%
\begin{eqnarray*}
\frac{{1}}{\sqrt{n}}\Bigl\Vert \max_{1\leq k\leq n}|S_{k}-M_{k}|
\Bigr\Vert_{2}&\leq&
\frac{{1}}{\sqrt{n}}\Bigl\Vert \max_{1\leq k\leq n}|M_{k}^{m}-M_{k}|
\Bigr\Vert _{2} \\
&&{}+\frac{{1}}{\sqrt{n}}\Vert \theta_{0}^{m}\Vert _{2}+\frac{{1}}{\sqrt{n}}\Bigl\Vert %
\max_{1\leq k\leq n}|\theta_{k}^{m}| \Bigr\Vert_{2}+\frac{{1}}{\sqrt{n}}\Bigl\Vert
\max_{1\leq k\leq n}|\overline{R}_{k}^{ m}| \Bigr\Vert_{2}.
\end{eqnarray*}
By Doob's maximal inequality for martingales and by stationarity, we conclude
that
\[
\frac{{1}}{\sqrt{n}}\Bigl\Vert \max_{1\leq k\leq n}|M_{k}^{m}-M_{k}|
\Bigr\Vert _{2}\leq
\Vert D_{0}^{m}-D_{0}\Vert _{2}.
\]
For $m$ fixed, since $(\theta_{k}^{m})_{k\in Z}$ is a stationary sequence
of square-integrable random variables, for any $A>0$, we have%
\begin{eqnarray*}
\frac{{1}}{n}\mathbb{E}\Bigl[\max_{1\leq k\leq n}|\theta_{k}^{m}|^{2}\Bigr]
&\leq&%
\frac{A^{2}}{n}+\frac{1}{n}\sum_{k=1}^{n}\mathbb{E[}|\theta
_{k}^{m}|^{2}I(|\theta_{k}^{m}|>A)] \\
&=&\frac{A^{2}}{n}+\mathbb{E[}|\theta_{0}^{m}|^{2}I(|\theta_{0}^{m}|>A)]
\end{eqnarray*}
and then, clearly,%
%
%e16 ###
\begin{equation}
\lim_{n\rightarrow\infty}\frac{{1}}{n}\mathbb{E}\Bigl[\max_{1\leq k\leq
n}|\theta_{k}^{m}|^{2}\Bigr]=0. \label{maxnegl}
\end{equation}
Then, taking into account \eqref{rests}, we easily obtain%
\[
\lim\sup_{n\rightarrow\infty}\frac{{1}}{\sqrt{n}}\Bigl\Vert \max_{1\leq
k\leq
n}|S_{k}-M_{k}| \Bigr\Vert_{2}\leq\Vert D_{0}^{m}-D_{0}\Vert _{2}+\Vert Y_{0}^{m}\Vert _{M^{+}}
\]
and the result follows by letting $m\rightarrow\infty$, from the fact
that $D_{0}^{m}\rightarrow D_{0}$ in $\mathbb{L}_{2}$. It is easy to see
that the martingale is unique.

\textit{Step \textup{3.} Necessity}.

Assume that the martingale approximation (\ref{maxcond}) holds. With the
notation $R_{n}=S_{n}-M_{n},$ we then have
\[
\lim_{n\rightarrow\infty}\frac{{1}}{\sqrt{n}}\Bigl\Vert \max_{1\leq k\leq
n}|R_{k}| \Bigr\Vert _{2}=0.
\]
In particular, this approximation implies that
%
%e17 ###
\begin{equation}
\lim_{n\rightarrow\infty}\frac{1}{\sqrt{n}}\max_{1\leq k\leq
n}\Vert \mathbb E(S_{k}|%
\mathcal{F}_{0})\Vert _{2}=0. \label{natural}
\end{equation}
From
\[
\Vert \overline{R}_{n}^{ n}\Vert _{2}\leq\Vert \mathbb
E(S_{n}|\mathcal{F}_{0})\Vert _{2},
\]
we deduce
that
\[
\Vert R_{n}^{n}\Vert _{2}=\Vert \theta_{0}^{n}-\theta_{n}^{n}+\overline{R}%
_{n}^{ n}\Vert _{2}\leq2\Vert \theta_{0}^{n}\Vert _{2}+\Vert \overline{R}%
_{n}^{ n}\Vert _{2}\leq3\max_{1\leq k\leq n}\Vert \mathbb E(S_{k}|\mathcal
{F}_{0})\Vert _{2},
\]
whence, by (\ref{natural}), it follows that
\[
\lim_{n\rightarrow\infty}\frac{\Vert R_{n}^{n}\Vert _{2}}{\sqrt{n}}=0.
\]

As a consequence, we obtain%
\[
\mathbb{E}(D_{0}^{n}-D_{0})^{2}=\frac{\mathbb
{E}(M_{n}^{n}-M_{n})^{2}}{n}=%
\frac{\mathbb{E}(R_{n}^{n}-R_{n})^{2}}{n}\rightarrow0\qquad\mbox{as }
n\rightarrow\infty.
\]
This shows that $D_{0}^{n}\rightarrow D_{0}$ in $\mathbb{L}_{2}.$
By the triangle inequality, followed by Doob's inequality, for any
positive integer $m$, we have
\begin{eqnarray*}
\frac{{1}}{\sqrt{n}}\Bigl\Vert \max_{1\leq k\leq n}|R_{k}^{m}| \Bigr\Vert _{2}
&\leq&
\frac{{1}}{\sqrt{n}}\Bigl\Vert \max_{1\leq k\leq n}|R_{k}| \Bigr\Vert _{2}+\frac
{{1}}{%
\sqrt{n}}\Bigl\Vert \max_{1\leq k\leq n}|M_{k}^{m}-M_{k}| \Bigr\Vert_{2} \\
&\leq&\frac{{1}}{\sqrt{n}}\Bigl\Vert \max_{1\leq k\leq n}|R_{k}| \Bigr\Vert %
_{2}+\Vert D_{0}^{m}-D_{0}\Vert .
\end{eqnarray*}
Now, letting $n\rightarrow\infty$ followed by $m\rightarrow\infty,$ we
obtain
%
%e18 ###
\begin{equation}
\lim_{m\rightarrow\infty}\lim\sup_{n\rightarrow\infty}\frac
{{1}}{\sqrt{n%
}}\Bigl\Vert \max_{1\leq k\leq n}|R_{k}^{m}| \Bigr\Vert _{2}=0. \label{rest1}
\end{equation}
Now, observe that by \eqref{rest}, $R_{n}^{m}-\overline{R}_{n}^{
m}=\theta
_{0}^{m}-\theta_{n}^{m}$. Then, for every fixed $m,$ by (\ref
{maxnegl}), we have
\[
\frac{{1}}{\sqrt{n}}\Bigl\Vert \max_{1\leq k\leq n}|\theta_{0}^{m}-\theta
_{k}^{m}| \Bigr\Vert _{2}\mathop{\rightarrow}_{n\rightarrow\infty}0.
\]
Thus, we conclude from \eqref{rest1} that
\[
\lim_{m\rightarrow\infty}{\Vert }Y_{0}^{ m}{\Vert }_{M^{+}}=0
\]
and the necessity follows.

%s3 ###
\section{Applications}

%s3.1 ###
\subsection{Applications using projective criteria}

The first application involves the class of variables satisfying the
Maxwell--Woodroofe condition~\cite{mw}.

\begin{proposition}
\label{underMW} Assume that
%
%e19 ###
\begin{equation}
\Delta(X_{0})=\sum_{k=1}^{\infty}\frac{\Vert \mathbb
{E}_{0}(S_{k})\Vert _{2}}{k^{3/2}%
}<\infty. \label{MW}
\end{equation}
The martingale approximation \eqref{maxcond} then holds.
\end{proposition}

\begin{pf}In order to verify condition (\ref{MAX}) of Theorem \ref{T},
we apply inequality (\ref{PU}) to the stationary sequence
$(Y_{k}^{m})_{k\in
\mathbb{Z}}$ defined by (\ref{defY}). Then
\[
\Biggl\Vert \max_{1\leq j\leq n}\Biggl|\sum_{k=0}^{j-1}Y_{k}^{m}\Biggr|
\Biggr\Vert_{2}\leq n^{1/2}\bigl(2\Vert Y_{0}^{m}\Vert _{2}+80\Delta(Y_{0}^{m})\bigr).
\]
First, note that by Peligrad and Utev \cite{PU1}, Proposition 2.5, we know
that condition (\ref{MW}) implies that $\Vert Y_{0}^{m}\Vert _{2}\rightarrow
0.$ We
complete the proof by showing that
\[
\Delta(Y_{0}^{m})\mathop{\longrightarrow}_{m\rightarrow\infty}0.
\]
Since $\Vert Y_{0}^{m}\Vert _{2}\rightarrow0,$ by the triangle inequality and
stationarity, every term of the series on the right-hand side of the
equality
\[
\Delta(Y_{0}^{m})=\sum_{k=1}^{\infty}\frac{1}{k^{3/2}}\Vert \mathbb{E}%
_{0}(Y_{0}^{m}+\cdots+Y_{k-1}^{m})\Vert _{2}
\]
tends to $0$ as $m\rightarrow\infty.$ Furthermore, because
\begin{eqnarray*}
\Vert \mathbb{E}_{0}(Y_{0}^{m}+\cdots +Y_{k-1}^{m})\Vert _{2}&=&\Biggl\Vert\mathbb{E}%
_{0}\Biggl(\frac{1}{m}\sum_{l=1}^{m}\sum_{i=0}^{k-1}\mathbb
{E}_{i}(X_{i+l})%
\Biggr)\Biggr\Vert _{2} \\
&\leq&\Vert  \mathbb{E}_{0}(X_{0}+\cdots+X_{k-1})\Vert _{2},
\end{eqnarray*}
each term in $\Delta(Y_{0}^{m})$ is dominated by the corresponding
term in $%
\Delta(X_{0})$, the latter being independent of $m$. The result
follows from the above considerations, along with the Lebesgue
dominated convergence theorem for the counting
measure.
\end{pf}

For the sake of applications, we give the following corollary.

\begin{corollary}
Assume that%
%
%e20 ###
\begin{equation}
\sum_{n=1}^{\infty}\frac{1}{\sqrt{n}}\Vert\mathbb{E}_{0}(X_{n})\Vert
_{2}<\infty. \label{mixg1}
\end{equation}
The martingale representation \eqref{maxcond} then holds.
\end{corollary}

The fact that (\ref{mixg1}) implies (\ref{MW}) was observed in Maxwell and
Woodroofe \cite{mw}.

We shall now combine Theorem \ref{T} with Rio's maximal inequality (\ref
{Rio}%
) to obtain the following proposition.

\begin{proposition}
\label{underRio}Assume that for any $j\geq0$,%
%
%e21 ###
\begin{equation}
\Gamma_{j}=\sum_{k\geq j}\Vert X_{j}\mathbb{E}_{0}(X_{k})\Vert _{1}<\infty\quad \mbox{and}\quad \frac{1}{m}\sum_{j=0}^{m-1}\Gamma_{j}\rightarrow0
\qquad\mbox{as }%
 m \rightarrow\infty. \label{condRio}
\end{equation}
The martingale representation \eqref{maxcond} then holds.
\end{proposition}

\begin{pf}
In order to verify condition (\ref{MAX}), we now apply the maximal
inequality (\ref{Rio}) to $(Y_{k}^{m})_{k\geq1}$ defined by (\ref{defY}).
We conclude that for $n\geq m$,
\begin{eqnarray*}
\Biggl\Vert \max_{1\leq j\leq n}\Biggl|\sum_{k=0}^{j-1}Y_{k}^{m}\Biggr|\Biggr \Vert
_{2}^{2}&\leq&
8n\Vert Y_{0}^{m}\Vert _{2}^{2}+16\sum_{j=1}^{n-1}\Vert Y_{0}^{m}\mathbb{E}%
_{0}(Y_{1}^{m}+\cdots+Y_{j}^{m})\Vert _{1} \\
&\leq&8n(12m+1)\Vert Y_{0}^{m}\Vert _{2}^{2}+16\sum
_{j=m+1}^{n-1}\Vert Y_{0}^{m}\mathbb{E}%
_{0}(Y_{m+1}^{m}+\cdots+Y_{j}^{m})\Vert _{1},
\end{eqnarray*}
where, in the last sum, we have implemented a decomposition into two
terms to deal
with overlapping blocks. So, for an absolute constant $C$,
\[
\frac{1}{n}\Biggl\Vert \max_{1\leq j\leq n}\Biggl|\sum_{k=0}^{j-1}Y_{k}^{m}\Biggr| \Biggr\Vert %
_{2}^{2}\leq C\Biggl(\frac{\Vert \mathbb{E}_{0}(S_{m})\Vert _{2}^{2}}{m}+\frac
{1}{n}%
\sum_{l=m+1}^{n-1}\Vert Y_{0}^{m}\mathbb{E}_{0}(Y_{m+1}^{m}+\cdots
+Y_{l}^{m})\Vert _{1}%
\Biggr).
\]
Since, for any $l>m$,
\begin{eqnarray*}
\Vert Y_{0}^{m}\mathbb{E}_{0}(Y_{m+1}^{m}+\cdots+Y_{l}^{m})\Vert _{1}& \leq&\frac
{1}{%
m}\sum_{j=1}^{m}\sup_{i>m}\Vert (\mathbb{E}_{0}(X_{j}))\mathbb{E}%
_{0}(X_{i}+\cdots+X_{i+l})\Vert _{1} \\
&\leq&\frac{1}{m}\sum_{j=1}^{m}\sum_{k\geq m}\Vert \mathbb
{E}_{0}(X_{j})\mathbb{E%
}_{0}(X_{k})\Vert _{1}
\end{eqnarray*}
and also
\[
\Vert \mathbb{E}_{0}(S_{m})\Vert _{2}^{2}\leq
2\sum _{j=0}^{m-1}\sum _{k=j}^{m-1}\Vert \mathbb{E}_{0}(X_{j})%
\mathbb{E}_{0}(X_{k})\Vert _{1},
\]
we then obtain, by the properties of
conditional expectations, that for a certain absolute constant~$C^{\prime}$,
\[
\frac{1}{n}\Biggl\Vert \max_{1\leq j\leq n}\Biggl|\sum
_{k=0}^{j-1}Y_{k}^{m}%
\Biggr| \Biggr\Vert _{2}^{2} \leq \frac{C^{\prime}}{m}\sum
_{j=0}^{m}\sum_{k\geq j}\Vert X_{j}\mathbb{E}%
_{0}(X_{k})\Vert _{1}
\]
and the result follows from condition (\ref{condRio}), by first letting
$%
n\rightarrow\infty$, followed by $m\rightarrow\infty$.
\end{pf}

The projective criteria in the next proposition were studied in
\cite{hh,HA,GG}, among others.

\begin{proposition}
Assume
%
%e22 ###
\begin{equation}
\mathbb E(X_{0}|\mathcal{F}_{-\infty})=0\quad  \mbox{almost surely
and}\quad
\sum_{i=1}^{\infty}\Vert\mathbb{E}_{-i}(X_{0})-\mathbb{E}%
_{-i-1}(X_{0})\Vert_{2}<\infty. \label{P}
\end{equation}
The martingale approximation \eqref{maxcond} then holds.
\end{proposition}

\begin{pf} The validity of this proposition easily follows by verifying
condition (\ref{MAX}) via maximal inequality (\ref{DM}) applied to $%
(Y_{k}^{m})_{k\geq1}$ defined by (\ref{defY}). Indeed, by (\ref{DM}),
the triangle inequality and stationarity, we have
\begin{eqnarray*}
\frac{1}{\sqrt{n}}\Biggl\Vert \max_{1\leq j\leq n}\Biggl|\sum_{k=0}^{j-1}Y_{k}^{m}\Biggr|\Biggr\Vert%
_{2}&\leq&2\sum_{i=0}^{\infty}\Vert\mathbb{E}_{-i}(Y_{0}^{m})-\mathbb
{E}%
_{-i-1}(Y_{0}^{m})\Vert_{2} \\
&\leq&\frac{2}{m}\sum_{i=0}^{\infty}\sum_{k=1}^{m}\Vert\mathbb
{E}_{-i}(X_{k})-%
\mathbb{E}_{-i-1}(X_{k})\Vert_{2}.
\end{eqnarray*}
Now, by stationarity, change of order of summation and change of
variable, %
\[
\frac{1}{\sqrt{n}}\Biggl\Vert\max_{1\leq j\leq n}\Biggl|\sum_{k=0}^{j-1}Y_{k}^{m}\Biggr|
\Biggr\Vert%
_{2}\leq\frac{2}{m}\sum_{k=1}^{m}\sum_{j=k}^{\infty}\Vert\mathbb{E}%
_{-j}(X_{0})-\mathbb{E}_{-j-1}(X_{0})\Vert_{2}.
\]
To verify condition (\ref{MAX}), we let $n\rightarrow\infty$ followed
by $%
m\rightarrow\infty.$ Note that the term on the right-hand side of the
previous inequality tends to $0$ as $m\rightarrow\infty$, by (\ref{P}).
\end{pf}

%s3.2 ###
\subsection{Application to mixing sequences}

The results in the previous section can be immediately applied to mixing
sequences, leading to the sharpest possible results and providing additional
information about the structures of these processes. Examples include
various classes of Markov chains and Gaussian processes.

We shall introduce the following mixing coefficients: for any two $%
\sigma$-algebras $\mathcal{A}$ and $\mathcal{B}$, define the strong mixing
coefficient $\alpha(\mathcal{A},\mathcal{B)}$,
\[
\alpha(\mathcal{A},\mathcal{B)=}\sup\{|\mathbb{P}(A\cap B)-\mathbb{P}(A)
\mathbb{P}(B)|;A\in\mathcal{A},B\in\mathcal{B\}},
\]
and the $\rho$-mixing coefficient, known also as the maximal
coefficient of correlation $\rho(\mathcal{A},\mathcal{B})$,
\[
\rho(\mathcal{A},\mathcal{B})=\sup\{\operatorname{Cov}(X,Y)/\Vert X\Vert
_{2}\Vert
Y\Vert_{2}\dvtx  X\in\mathbb{L}_{2}(\mathcal{A}),Y\in\mathbb{L}_{2}(
\mathcal{B})\}.
\]

For the stationary sequence of random variables $(X_{k})_{k\in
\mathbb{Z}},$ we also define $\mathcal{F}_{m}^{n}$, the $\sigma$-field generated
by $X_{i}$ with indices $m\leq i\leq n$. $\mathcal{F}^{n}$ denotes the $
\sigma$-field generated by $X_{i}$ with indices $i\geq n$ and $\mathcal
{F}%
_{m}$ denotes the $\sigma$-field generated by $X_{i}$ with indices
$i\leq m.
$ The sequences of coefficients~$\alpha(n)$ and $\rho(n)$ are then defined
by
\[
\alpha(n)=\alpha(\mathcal{F}_{0},\mathcal{F}_{n}^{n})\quad \mbox{and}\quad %
\rho(n)=\rho(\mathcal{F}_{0},\mathcal{F}^{n}).
\]

Equivalently (see \cite{rick}, Chapter 4),
\[
\rho(n)=\sup\{\Vert\mathbb{E}(Y|\mathcal{F}_{0})\Vert_{2}/\Vert Y\Vert
_{2}\dvtx  Y\in\mathbb{L}_{2}(\mathcal{F}^{n}),\mbox{ }\mathbb{E}(Y)=0\}.
\]

Finally, we say that the stationary sequence is strongly
mixing if $\alpha(n)\rightarrow0$ as $n\rightarrow\infty$ and $\rho$-mixing
if $\rho(n)\rightarrow0$ as $n\rightarrow\infty$.

An interesting application of Proposition \ref{underMW} is to
$\rho$-mixing sequences. It is well known that the central limit
theorem and its invariance principle hold for stationary centered
sequences with finite second moments under the
assumption
%
%e23 ###
\begin{equation}
\sum_{k=1}^{\infty}\rho(2^{k})<\infty, \label{condrho}
\end{equation}
where $\rho(n)=\rho(\mathcal{F}_{0},\mathcal{F}^{n}\mathcal{)}$.
Let us
recall that the central limit theorem is due to \cite{II}, while
the invariance principle is found in \cite{pm,shao,u89,u91}. The fact that condition (\ref{condrho}) is sharp in this
context is due to \cite{rick}, Volume 1, page 367, and Volume
3, Theorem 34.13.
Bradley's example shows that if (\ref{condrho}) fails, then $%
S_{n}/ \Vert S_{n}\Vert _2$ might have
non-degenerate non-normal distributions as weak limit points.

As a corollary of Proposition \ref{underMW}, we obtain the conditional
invariance principle for $\rho$-mixing sequences.

\begin{proposition}
Assume $\sum_{k=1}^{\infty}\rho(2^{k})<\infty.$ The martingale
representation \eqref{maxcond} then holds.
\end{proposition}

\begin{pf} As in \cite{mpu2},
for a positive
constant $C$, we have
\[
\sum_{r=0}^{\infty}\frac{\Vert
\mathbb{E}(S_{2^{r}}|\mathcal{F}_{0})\Vert_{2}}{2^{r/2}}\leq C
\Vert X_0\Vert _2\sum_{j=0}^{\infty}\rho(2^{j}).
\]
\upqed\end{pf}

To obtain sharp results for strongly mixing sequences, we shall
use Proposition \ref{underRio}.

According to Doukhan, Massart and Rio \cite{DMR}, a condition that is optimal
for CLT or the invariance principle for strongly mixing sequences is
%
%e24 ###
\begin{equation} \label{condalha}
\sum _{k\geq1}{\mathbb{E}}X_{0}^{2}I\bigl(|X_{0}|\geq
Q_{|X_{0}|}(2\alpha_{k})\bigr)<\infty,
\end{equation}
where $Q_{|X_{0}|}$ denotes the cadlag inverse
of the function $t\rightarrow P(|X_{0}|>t).$ Also under this
condition, we add the additional information given by Theorem
\ref{T}.

\begin{proposition}
Assume that condition \eqref{condalha} is satisfied. The martingale
representation \eqref{maxcond} then holds.
\end{proposition}

\begin{pf}
We shall just verify the condition of Proposition \ref{underRio}. Note that
on the set $[0,{\ P}(|Y|>0)]$, the function $H_{Y}\dvtx x\rightarrow
\int_{0}^{x}Q_{Y}(u)\,\mathrm{d}u$ is an absolutely continuous and increasing function
with values in $[0,{\ E}|Y|]$. Denote by $G_{Y}$ the inverse of
$H_{Y}$. With this notation, by Merlev\`{e}de and Peligrad \cite{mp},
relation (4.84), we have
\[
\Vert X_{j}\mathbb{E}(X_{k}|\mathcal{F}_{0})\Vert _{1}\leq3\int_{0}^{\Vert
{\mathbb{%
E}}(X_{k}|\mathcal{F}_{0})\Vert_{1}}Q_{|X_{0}|}\circ G(u)\,\mathrm{d}u
\]
and we then majorize the right-hand side in the previous inequality by
Dedecker and Doukhan~\cite{dd}, Proposition 1, to obtain
\[
\Vert X_{j}\mathbb{E}(X_{k}|\mathcal{F}_{0})\Vert _{1}\leq6\int_{0}^{2\alpha
(k)}Q_{|X_{0}|}^{2}\,\mathrm{d}u .
\]
Therefore,%
\begin{eqnarray*}
\sum_{k\geq j}\Vert X_{j}\mathbb{E}_{0}(X_{k})\Vert _{1}&\leq&6\sum_{k\geq
j}\int_{0}^{2\alpha(k)}Q_{|X_{0}|}^{2}\,\mathrm{d}u  \\
&\leq&6\sum_{k\geq j}{\mathbb{E}}X_{0}^{2}I\bigl(|X_{0}|\geq Q_{|X_{0}|}(2\alpha
_{k})\bigr) \rightarrow0\qquad \mbox{as }j\rightarrow\infty
 .
\end{eqnarray*}
\upqed\end{pf}

Note that the coefficient $\alpha(k)$ is defined by using only one
variable in the future. Moreover, by the Cauchy--Schwarz inequality,
condition (%
\ref{condalha}) is satisfied if the variables have finite moments of
order $%
2+\delta$ for a $\delta>0$ and
\[
\sum_{k\geq1}\alpha(k)^{\delta/(2+\delta)}<\infty.
\]
An excellent source of information for classes of mixing sequences and
classes of Markov chains satisfying mixing conditions is the book by
Bradley %
\cite{rick}. Further applications can be obtained by using the coupling
coefficients in \cite{dp}.

%s3.3 ###
\subsection{Application to additive functionals of reversible Markov
chains}

For reversible Markov processes (i.e., $Q=Q^{\ast}$), the invariance principle
under an optimal condition is known since Kipnis and Varadhan \cite
{KV}. The following
is a formulation in terms of martingale approximation.

\begin{proposition}
Let $(\xi_{i})_{i\in\mathbb{Z}}$ be a stationary reversible
Markov chain and $f\in\mathbb{L}_{0}^{2}(\pi)$ with the property
%
%e25 ###
\begin{equation}
\lim_{n\rightarrow\infty}\frac{\operatorname{var}(S_{n})}{n}\rightarrow\sigma
_{f}^{2}<\infty%
.
\label{revcondCLT}
\end{equation}
The martingale approximation satisfying \eqref{maxcond}
then holds.
\end{proposition}

\begin{pf}
We have to verify condition \eqref{MAX}. Denote by $\rho_{f}$ the
spectral measure of $f$ corresponding to the self-adjoint operator
$Q$ on $\mathbb{L}_{2}(\pi)$. It is well known that the
assumption \eqref{revcondCLT} for $f \in\mathbb L^2_0$ implies
that $ \int _{-1}^{1}(1-t)^{-1}\rho_{f}(\mathrm{d}t)<\infty$ (see
\cite{KV}). Define $Y_{0}^{m}$ by \eqref{defY}.
By the maximal inequality~\eqref{LW}, we have%
\[
\frac{1}{n}\mathbb{E}\Biggl(\max_{1\leq j\leq n}\Biggl|\sum
_{k=0}^{j-1}Y_{k}^{m}\Biggr|\Biggr)%
^{2}\leq27\sum_{k\geq0}\mathbb{E(}Y_{0}^{m}Y_{k}^{m}),
\]
provided that the sum on the right-hand side is finite.
To prove it, by using spectral calculus for
the self-adjoint operator $Q$, %(see Kipnis and Varadhan, 1986)
we obtain
\[
\sum_{k\geq0}\mathbb{E(}Y_{0}^{m}Y_{k}^{m})\leq\frac{1}{m^{2}}%
\int _{-1}^{1}\frac{(1+t+\cdots+t^{m-1})^{2}}{(1-t)}\rho_{f}(\mathrm{d}t)
\]
and therefore, for every positive integer $m>0,$%
\[
\Vert Y_{0}^{m}\Vert _{M^{+}}^{2}\leq27\int _{-1}^{1}\frac{%
(1+t+\cdots+t^{m-1})^{2}}{m^{2}(1-t)}\rho_{f}(\mathrm{d}t).
\]
Since $%
\int _{-1}^{1}(1-t)^{-1}\rho_{f}(\mathrm{d}t)<\infty$, the right-hand
side is finite and, by the dominated convergence theorem,%
\[
\lim_{m\rightarrow\infty}\Vert Y_{0}^{m}\Vert _{M^{+}}^{2}=0.
\]
\upqed\end{pf}

Similar results are expected to hold for other classes of stationary and
ergodic Markov chains when $Q~$is not necessarily self-adjoint, but
instead satisfies a quasi-symmetry or strong sector condition, or is
symmetrized.
See \cite{wu} and \cite{svy} for
these related
processes.
%s3.4 ###
\subsection{Application to additive functionals of normal Markov chains}

For additive functionals of normal Markov chains $(QQ^{\ast}=Q^{\ast
}Q),$ the
central limit theorem below is a result of Gordin and Lifshitz \cite{gl}.
As an application of Theorem \ref{T2}, we give an alternative proof.

Let $\rho_{f}$ be the spectral measure on the closed unit disk $D
\subset\mathbb C$ corresponding to the function
$f\in\mathbb{L}_{0}^{2}(\pi)$.

\begin{proposition}
\label{normal} Let $(\xi_{i})_{i\in\mathbb{Z}}$ be a stationary normal
Markov chain and
a function $f\in\mathbb{L}_{0}^{2}(\pi),$
satisfying the condition
%
%e26 ###
\begin{equation}
\int _{D}\frac{1}{|1-z|}\rho_{f}(\mathrm{d}z)<\infty.
\label{normcond}
\end{equation}
The martingale approximation \eqref{MA} then holds.
\end{proposition}

\begin{pf}
According to
Theorem \ref{T2}, we have to verify condition (\ref{MA1}). By using
spectral calculus as in \cite{1994}, Chapter 4, after some
computations, we get%
\[
\lim\sup_{n\rightarrow\infty}\frac{1}{n}\Biggl\Vert %
\sum_{k=0}^{n-1}Y_{k}^{m}\Biggr\Vert _{2}^{2}\leq4\int _{D}\frac{%
|1+z+\cdots+z^{m-1}|^{2}}{m^{2}|1-z|}\rho_{f}(\mathrm{d}z)
\]
and condition (\ref{MA1}) is therefore satisfied by condition (\ref
{normcond}) and the dominated convergence theorem.
\end{pf}

Condition \eqref{normcond} has an interesting equivalent
formulation in terms of conditional moments that is in the
spirit of (and which implies) the Mawxell--Woodroofe
condition \eqref{MW}.

\begin{remark}
Condition \eqref{normcond} is equivalent to
%
%e27 ###
\begin{equation}
\sum_{k=1}^{\infty}\frac{\Vert \mathbb
{E}_{0}(S_{k})\Vert _{2}^{2}}{k^{2}}<\infty%
. \label{Gap}
\end{equation}

Condition \eqref{Gap} is further implied by
%
%e28 ###
\begin{equation}
\sum_{k=1}^{\infty}\Vert \mathbb{E}_{0}(X_{k})\Vert _{2}^{2}<\infty.
\label{cor2}
\end{equation}
\end{remark}

The equivalence in the above remark can be found in \cite{Cuny}, Lemma 2.1.
The fact that \eqref{cor2} implies~\eqref{Gap} is easily established,
much like
the proof that \eqref{mixg1} implies \eqref{MW}.

%s4 ###
\section*{Acknowledgements}

Mikhail Gordin was supported in part by a Charles Phelps Taft
Memorial Fund grant and RFBR Grant 10-01-00242\_a. Magda Peligrad
was supported in part
by a Charles Phelps Taft Memorial Fund grant and NSA Grant
H98230-09-1-0005.
The authors are grateful to the referees for carefully reading the
paper and
for numerous suggestions that improved the presentation of the paper.

\printhistory

\end{document}